\documentclass[12pt,leqno]{amsart}

\usepackage{amsmath, amssymb,latexsym}
\usepackage{amsthm,amscd,amsfonts}
\usepackage{graphicx} 

\begin{document}

\newtheorem{thm}{Theorem}[section]
\newtheorem{prop}[thm]{Proposition}
\newtheorem{lem}[thm]{Lemma}
\newtheorem{cor}[thm]{Corollary}
\newtheorem{conj}[thm]{Conjecture}
\newtheorem{ddef}{Definition}
\newtheorem{ex}[thm]{Example}
\newtheorem{rem}[thm]{Remark}
\newtheorem{notation}[thm]{Notation}
\numberwithin{equation}{section}
\newcommand{\bthm}{\begin{thm}} \newcommand{\ethm}{\end{thm}}
\newcommand{\bthms}{\begin{thm*}} \newcommand{\ethms}{\end{thm*}}
\newcommand{\blem}{\begin{lem}} \newcommand{\elem}{\end{lem}}
\newcommand{\bcor}{\begin{cor}} \newcommand{\ecor}{\end{cor}}
\newcommand{\bprop}{\begin{prop}} \newcommand{\eprop}{\end{prop}}
\newcommand{\bproof}{\begin{proof}} \newcommand{\eproof}{\end{proof}}
\newcommand{\bddef}{\begin{ddef}} \newcommand{\eddef}{\end{ddef}}
\newcommand{\bconj}{\begin{conj}} \newcommand{\econj}{\end{conj}}
\newcommand{\brem}{\begin{rem}} \newcommand{\erem}{\end{rem}}
\newcommand{\bca}{\begin{cases}} \newcommand{\eca}{\end{cases}}
\newcommand{\beq}{\begin{equation}} \newcommand{\eeq}{\end{equation}}
\newcommand{\beqs}{\begin{equation*}} \newcommand{\eeqs}{\end{equation*}}
\newcommand{\beqa}{\begin{eqnarray}} \newcommand{\eeqa}{\end{eqnarray}}
\newcommand{\beqas}{\begin{eqnarray*}} \newcommand{\eeqas}{\end{eqnarray*}}
\newcommand{\barr}{\begin{array}} \newcommand{\earr}{\end{array}}
\newcommand{\btab}{\begin{tabular}} \newcommand{\etab}{\end{tabular}}
\newcommand{\bit}{\begin{itemize}} \newcommand{\eit}{\end{itemize}}
\newcommand{\ben}{\begin{enumerate}} \newcommand{\een}{\end{enumerate}}
\newcommand{\bce}{\begin{center}} \newcommand{\ece}{\end{center}}
\newcommand{\defeq}{\stackrel{\rm def}{=}}
\newcommand{\bd}{\partial}
\newcommand{\opp}{\mathrm{opp}}
\newcommand{\hz}{\widehat{0}}
\newcommand{\ho}{\widehat{1}}
\newcommand{\wh}{\widehat}
\newcommand{\wt}{\widetilde}
\newcommand{\cov}{\lhd}
\newcommand{\cocov}{\rhd}
\renewcommand{\AA}{\mathcal{A}}
\newcommand{\BB}{\mathcal{B}}
\newcommand{\CC}{\mathcal{C}}
\newcommand{\DD}{\mathcal{D}}
\newcommand{\EE}{\mathcal{E}}
\newcommand{\FF}{\mathcal{F}}
\newcommand{\II}{\mathcal{I}}
\newcommand{\LL}{\mathcal{L}}
\newcommand{\NN}{\mathcal{N}}
\newcommand{\HH}{\mathcal{H}}
\newcommand{\RR}{\mathcal{R}}
\newcommand{\MM}{\mathcal{M}}
\newcommand{\QQ}{\mathcal{Q}}
\renewcommand{\SS}{\mathcal{S}}
\newcommand{\al}{\alpha}
\newcommand{\be}{\beta}
\newcommand{\de}{\delta} \newcommand{\De}{\Delta} 
\newcommand{\ga}{\gamma} \newcommand{\Ga}{\Gamma}
\newcommand{\la}{\lambda} \newcommand{\La}{\Lambda}
\newcommand{\om}{\omega} \newcommand{\Om}{\Omega}
\newcommand{\si}{\sigma} \newcommand{\Si}{\Sigma}
\newcommand{\ze}{\zeta}
\newcommand{\vphi}{\varphi}
\newcommand{\eps}{\varepsilon}
\def\N{\mathbb{N}}
\def\P{\mathbb{P}}
\def\Z{\mathbb{Z}}
\newcommand{\Q}{\mathbb{Q}}
\newcommand{\R}{\mathbb{R}}
\newcommand{\C}{\mathbb{C}}
\newcommand{\F}{\mathbb{F}}
\renewcommand{\L}{\mathbb{L}}
\def\bb{\mathbf{b}}
\def\cc{\mathbf{c}}
\def\ff{\mathbf{f}}
\def\bg{\mathbf{g}}
\def\hh{\mathbf{h}}
\def\kk{\mathbf{k}}
\newcommand{\st}{\,:\,} 
\newcommand{\sbseq}{\subseteq}
\newcommand{\spseq}{\supseteq}
\newcommand{\larr}{\leftarrow}
\newcommand{\rarr}{\rightarrow}
\newcommand{\Larr}{\Leftarrow}
\newcommand{\Rarr}{\Rightarrow}
\newcommand{\lrarr}{\leftrightarrow}
\newcommand{\Lrarr}{\Leftrightarrow}
\def\rank{\text\rm{rank}}
\def\Lraw{\Longrightarrow}
\def\lraw{\longrightarrow}
\def\Llaw{\Longleftarrow}
\def\llaw{\longleftarrow}
\def\Llraw{\Longleftrightarrow}
\def\llraw{\longleftrightarrow}
\def\wtx{\underset{\displaystyle{\widetilde{}}}{x}}
\def\wth{\underset{\displaystyle{\widetilde{}}}{h_0}}
\def\({\left(}
\def\){\right)}
\def\no={\,{\,|\!\!\!\!\!=\,\,}}

\title[A $q$-analogue of the FKG inequality]
{A $q$-analogue of the FKG inequality \\ and some applications}
\author{Anders Bj\"orner}

\address{Institut Mittag-Leffler, Aurav\"agen 17, S-182 60 Djursholm, Sweden}
\email{bjorner@mittag-leffler.se}
\address{Kungl. Tekniska H\"ogskolan, Matematiska Inst.,
 S-100 44 Stockholm, Sweden}
\email{bjorner@math.kth.se}
\thanks{Research supported by the Knut and Alice Wallenberg Foundation,
grant KAW.2005.0098}

\begin{abstract} 
Let $L$ be a finite distributive lattice and $\mu : L \rightarrow {\mathbb R}^{+}$
a log-supermodular function. 
For functions $k: L \rightarrow {\mathbb R}^{+}$ let
$$E_{\mu} (k; q) \defeq \sum_{x\in L} k(x) \mu (x) q^{{\mathrm rank}(x)} \in {\mathbb R}^{+}[q].$$
We prove for any pair $g,h: L\rightarrow {\mathbb R}^{+}$ of monotonely increasing functions, that
$$E_{\mu} (g; q)\cdot E_{\mu} (h; q)\   \ll \ E_{\mu} (1; q)\cdot E_{\mu} (gh; q),
$$
where ``$ \, \ll \, $'' denotes coefficientwise inequality of real polynomials.
The FKG inequality of Fortuin, Kasteleyn and Ginibre (1971) is 
the real number inequality obtained by specializing  to $q=1$.

The polynomial FKG inequality has applications to $f$-vectors of joins 
and intersections of
simplicial complexes, to Betti numbers of intersections of certain Schubert
varieties, and to the following kind of
 correlation inequality for  power series weighted by Young tableaux.

Let $Y$ be the set of all integer partitions.
Given functions $k,  \mu: Y \rarr \R^+$, and 
parameters $0\le s\le t$, define the formal power series
$$F_{\mu}(k ; z) \defeq \sum_{\la\in Y} k(\la) \, \mu(\la)\,  (f_{\la})^t \, \frac{z^{|\la|}}{(|\la| !)^s} 
\, \in \R^+ [[z]], $$
where $f_{\la}$ is the number of standard Young tableaux of shape $\la$.
Assume that $\mu: Y\rarr \R^+$ is log-supermodular, and that
$g, h: Y \rarr \R^+$ are  monotonely increasing with respect to containment order of
partition shapes. Then 
$$F_{\mu}(g;z) \, \cdot \,  F_{\mu}(h;z) \, \ll \, F_{\mu}(1;z)  \, \cdot \, F_{\mu}(gh;z).
$$

\end{abstract}

\maketitle

\section{Introduction}

Suppose that $L$ is a finite lattice, and for functions $k, \mu: L \rarr \R$ let
\beqs E_{\mu}[k] \defeq 
\sum_{x\in L} k(x) \mu(x). 
\eeqs
A function $g: L \rarr \R$  is {\em increasing} 
if $x\le y$ implies $g(x)\le g(y)$, and 
{\em decreasing} 
if $x\le y$ implies $g(x)\ge g(y)$.
Two functions $g, h: L \rarr \R$  are {\em comonotone} if either both are increasing 
or else  both are decreasing. They are
{\em countermonotone} if one is increasing and the other decreasing.

Let $\R^+ \defeq \{r\in\R\st r\ge 0\}$.
A  function $\mu: L\rarr \R^+$ is said to be  {\em log-supermodular} if
$$\mu(x) \mu(y) \le \mu(x\wedge y)\mu(x\vee y), \mbox{ for all } x,y\in L. $$

The  following theorem is due to Fortuin, Kasteleyn and Ginibre 
\cite{FKG}. 

\bthm (FKG inequality) \label{fkg}Let $L$ be a finite distributive lattice,  \label{FKG}
$\mu: L\rarr \R^+$ a log-supermodular weight function, and $g$ and $h$ 
comonotone functions 
$g,h: L\rarr \R^+$. Then
$$E_{\mu} [g]\cdot E_{\mu} [h]\ \le\ E_{\mu} [1]\cdot E_{\mu}[gh].
$$
For countermonotone functions the inequality is reversed.
\ethm

The FKG inequality arose as a correlation inequality in statistical
mechanics, more precisely in the study of Ising ferromagnets
and the random cluster model.
It  has found many applications in extremal
and probabilistic combinatorics. See e.g. \cite {AS} for an expository
 account and 
 pointers to related results.
 
 In  this paper we prove a polynomial inequality that implies the  original FKG
 inequality via specialization, and which provides more
 detailed information. The basic version is given in 
 Section 2 and a more elaborate one in Section 5.
  
The polynomial FKG inequality has applications to $f$-vectors of joins 
and intersections of
simplicial complexes and to Betti numbers of intersections of certain Schubert
varieties. These applications appear in Sections 3 and 4.

In Section 6 we derive 
a class of correlation inequalities for  power series on partitions weighted by Young tableaux.
This class contains
the poissonization of Plancherel measure, a probability measure on the set 
 $Y$  of all integer partitions discussed in \cite{BOO}. 


\section{A $q$-FKG inequality}

Only the most basic facts about distributive  lattices are used in this
paper. Standard facts and notation can be found in e.g. 
 \cite{Bir} or \cite{EC1}. For instance, we use that a finite such  lattice has 
 bottom and top elements $\hz$ and $\ho$ and a {\em rank function}
 $r$, such that $r(x)$   equals the length of any maximal chain in the interval
 $[\hz, x]$. Furthermore, this function satisfies the {\em modular law}
 \beq\label{modular}
 r(x)+r(y)=r(x\wedge y)+r(x\vee y), \quad \mbox{ for all $x,y\in L$.}
 \eeq

For polynomials $p(q), s(q) \in \R[q]$,
$$\mbox{ let $p(q) \ll s(q)$
mean that $s(q)-p(q)\in \R^+ [q]$.}$$

Let $L$ be a finite distributive lattice  with rank function $r$. 
For functions $\mu, k: L\rarr \R^+$ define the polynomial
$$E_{\mu}(k; q)\defeqß\sum_{x\in L} k(x) \,  \mu(x)\, q^{r(x)} \in \R^+ [q] .$$

\bthm \label{qfkg}
Let $L$ be a finite distributive lattice, 
$\mu: L\rarr \R^+$ a log-supermodular weight function
and $g$ and $h$ 
comonotone functions 
$g,h: L\rarr \R^+$. Then
$$E_{\mu} (g; q)\cdot E_{\mu} (h; q)\   \ll \ E_{\mu} (1; q)\cdot E_{\mu} (gh; q).
$$
For countermonotone functions the inequality is reversed.
\ethm

The original FKG inequality (Theorem \ref{fkg}) is obtained by putting $q=1$.
Slightly  more general $q$-FKG inequalies are given in Theorem \ref{qfkg.alt}.
\bproof
We may assume that both $g$ and $h$ are increasing functions on $L$. The cases
when one or both are decreasing are easily deduced from this.

Let
$$ \Phi(q) \defeq \ E_{\mu} (1; q)\cdot E_{\mu} (gh; q) - E_{\mu} (g; q)\cdot E_{\mu} (h; q) ,
$$
and for $x,y\in L$ let
$$\phi(x,y)\defeq \mu(x)\mu(y)[g(x)-g(y)][h(x)-h(y)].
$$
Note that $\phi(x,y)=\phi(y,x)$.
A simple computation, consisting of grouping all the terms of $\Phi(q)$
that contain $\mu(x)\mu(y)$ together, shows that
$$\Phi(q)=\sum_{\{x,y\} \in \binom{L}{2}} \phi(x,y) \ q^{r(x)+r(y)}.
$$
Hence, the degree $d$ coefficient of $\Phi(q)$ equals
\beq\label{dcoeff}
 \Phi_d = \sum_{\{x,y\}\in \binom{L}{2}, r(x)+r(y)=d} \phi(x,y).
\eeq

Now consider another, slightly coarser, grouping of terms.
For $u,v\in L$ such that $u\le v$, let $C(u,v)$ denote the set of all pairs of
relative complements in the interval $[u,v]$, that is, unordered pairs
$\{x,y\}$ such that $x\wedge y=u, x\vee y=v$. Let
\beq\label{def2} \psi(u,v)\defeq  \sum_{\{x,y\}\in C(u,v)} \phi(x,y).
\eeq
It follows from equation (\ref{dcoeff}) and the modular law (\ref{modular}) for the
rank function  that 
\beq\label{eq7} 
\Phi_d =  \sum_{u\le v, r(u)+r(v)=d} \psi(u,v).
\eeq
Our aim is to prove the following.
\vspace{3mm}

\noindent
{\bf Claim.} {\em For all $u \le v$ in $L$}
\beq \label{claim}
\psi (u,v) \ge 0.
\eeq
The nonegativity of $\Phi_d$ then follows. 

To prove Claim (\ref{claim}) we may without loss of generality make the
following three assumptions:
\ben
\item[(i)] $\mu(x)> 0$ for all $x\in L$,
\item[(ii)] $u=\hz$ and $v=\ho$,
\item[(iii)] $L$ is Boolean.
\een
The reasons that these simplifications are valid are as follows.

(i) Let $L' \defeq \{x\in L\st \mu(x)>0 \}$. If $x\in L'$ and $y\in L'$ then log-super\-modularity
implies that also $x\wedge y \in L'$ and $x\vee y \in L'$. Hence, $L'$ is a
sublattice. Since elements outside $L'$ don't contribute to the sums involved
in the statement of the theorem we may pass from $L$ to consideration of the
distributive sublattice $L'$.

(ii) Every interval $[u,v]$ is itself a finite distributive lattice, and what
 goes into Claim (\ref{claim}) depends only on this substructure.

(iii) 
It is known that the set $C(\hz, \ho) $ of complemented elements in a finite distributive lattice $L$ 
forms  a Boolean sublattice, see \cite[p. 18]{Bir}.
By definition (\ref{def2}), $\psi (\hz, \ho) $ only depends on the restriction
to this Boolean sublattice. 
\medskip

After these simplifications we return to the 
proof of Claim (\ref{claim}), for which it remains to show that 
\beq\label{eq3}\psi (\hz, \ho) \ge 0
\eeq
in a finite Boolean lattice $L$. This will be done by induction on the rank
of $L$. The inequality is clearly valid for $\rank(L)=1$, in which case
$L=\{\hz, \ho\}$.

Suppose then
that $L$ is Boolean, and for each element $x\in L$, let $x^c$ denote
its complement in $L$.
Let $B\defeq [\hz, a^c]$ for some atom (rank one element) $a\in L$.

Let $\widetilde{\mu}(x)  \defeq   \mu(x) \mu(x^c)$. Note that 
$\wt{\mu}$ is a log-supermodular 
function on $L$, and hence in particular on $B$, and that
\beq\label{eq4}
\sum_{x\in B} \wt{\mu}(x) = \sum_{x\in B} \wt{\mu}(x^c)=\sum_{x\in B} \wt{\mu}(x\vee a) .
\eeq

Let
$$ \al \defeq \sum_{x\in B} [g(x)-g(x^c)] \,  \wt{\mu}(x) 
\mbox {\quad and \quad} \be \defeq \sum_{x\in B} [h(x)-h(x^c)] \,  \wt{\mu}(x) .
$$
The plan is to show that
\beq\label{eq5}
\psi (\hz, \ho) \ge \al\cdot\be,
\eeq
and
\beq\label{eq6}
\al\le 0  {\quad and \quad} \be\le 0,
\eeq
from which inequality (\ref{eq3}) follows. 

The interval $B$ is a Boolean lattice of 
rank one less than $L$. Hence, by the induction hypothesis 
$\psi(u,v)\ge 0$ for all $u\le v$ in $B$,  
which means that the statement of the theorem is valid for $B$.
This statement, specialized to $q=1$, shows that 
for any pair of increasing nonnegative functions $s, t$ on $B$:
\beqa\label{induct}
\sum_{x\in B} s(x) \,  \wt{\mu}(x)  \cdot
\sum_{x\in B} t(x) \,  \wt{\mu}(x)  \le
\sum_{x\in B}  \,  \wt{\mu}(x)  \cdot
\sum_{x\in B} s(x)t(x) \,  \wt{\mu}(x).
\eeqa


Now, define two functions $B\rarr \R^+$ this way:
\begin{eqnarray*}
\widetilde{g}(x) &\defeq &  e+g(x)-g(x^c) \\
\wt{h}(x) & \defeq &  e+h(x)-h(x^c),
\end{eqnarray*}
where $e$ is some number such that $e\ge  g(\ho), h(\ho)$. 
Then $\wt{g}$ and  $\wt{h}$
are nonnegative and increasing. Hence, by the induction hypothesis (\ref{induct})
\beqas\label{induct.1}
\sum_{x\in B} \wt{g}(x) \,  \wt{\mu}(x)  \cdot
\sum_{x\in B} \wt{h}(x) \,  \wt{\mu}(x)  \le
\sum_{x\in B}  \,  \wt{\mu}(x)  \cdot
\sum_{x\in B} \wt{g}(x)\wt{h}(x) \,  \wt{\mu}(x).
\eeqas
After expanding these products and cancelling terms this inequality 
reduces to
\beqas\label{eq1}
\sum_{x\in B} [g(x)-g(x^c)] \,  \wt{\mu}(x)  \cdot
\sum_{x\in B} [h(x)-h(x^c)]  \,  \wt{\mu}(x)  
\eeqas
$$ \le
\sum_{x\in B} [g(x)-g(x^c)][h(x)-h(x^c)]  \,  \wt{\mu}(x),
$$
which is inequality (\ref{eq5})


Next, let
$$\ga(x) \,\defeq \, \frac{\wt{\mu}(x\vee a)}{\wt{\mu}(x)}.
$$
Since ${\wt{\mu}(x)} > 0$
this defines a function $\ga: B\rarr \R^+$, which due to the 
log-supermodularity of $\wt{\mu}$ is increasing:
$$x<y \mbox{ and } y\in B \;\;\Rightarrow\;\; \wt{\mu}(y) \wt{\mu}(x\vee a) 
\le \wt{\mu}(x) \wt{\mu}(y\vee a) 
\;\;\Rightarrow\;\;  \ga(x) \le \ga(y).
$$
By the induction hypothesis (\ref{induct})
\beqas\label{induct.1}
\sum_{x\in B} g(x) \,  \wt{\mu}(x)  \cdot
\sum_{x\in B} \ga(x) \,  \wt{\mu}(x)  \le
\sum_{x\in B}  \,  \wt{\mu}(x)  \cdot
\sum_{x\in B} g(x) \ga(x) \,  \wt{\mu}(x),
\eeqas
which because of equation  (\ref{eq4}) reduces to 
\beqas
 \sum_{x\in B} g(x) \,  \wt{\mu}(x)  \le
\sum_{x\in B} g(x) \,  \wt{\mu}(x\vee a) .
\eeqas
This together with
\beqas
\sum_{x\in B} g(x) \,  \wt{\mu}(x\vee a)  \le
\sum_{x\in B} g(x\vee a) \,  \wt{\mu}(x\vee a) 
= \sum_{x\in B} g(x^c) \,  \wt{\mu}(x^c)
= \sum_{x\in B} g(x^c)  \wt{\mu}(x) 
\eeqas
proves inequality (\ref{eq6}) for $\al$. 
The analogous argument for $\be$ then finishes the proof.
\eproof

\section{Face numbers of intersections of simplicial complexes}
\noindent

Let $\De$ and $\Ga$ be simplicial complexes on the same finite vertex set $V$.
Define the 
{\em $f$-polynomial}  $f_{\De}(q)=\sum_{i\ge 0} f_{i-1}(\De)\, q^{i}$,
where $f_{i-1}(\De)$ denotes the number of $(i-1)$-dimensional faces,
and similarly for $\Ga$.

\bthm
$$ f_{\De}(q) \cdot f_{\Ga}(q)  \ll (1+q)^{|V|} \cdot f_{\De \cap \Ga}(q)
$$
\ethm
\bproof Apply Theorem \ref{qfkg} to the Boolean lattice $2^V$, putting
$$\mu(x)\equiv 1, \quad\quad 
g(x)=\begin{cases} 
1, \mbox{ if } x\in\De \\ 
0, \mbox{ if } x\not\in\De \end{cases} , \mbox{ and } \quad
  h(x)=\begin{cases} 1, \mbox{ if } x\in\Ga \\ 0, \mbox{ if } x\not\in\Ga \end{cases} $$
\eproof

The theorem can equivalently be formulated in terms of the simplicial join operation:
$$ f_i(\De \ast \Ga) \;\le\;  f_i(2^V  \ast (\De\cap\Ga)), \mbox{ for all $i$.}
$$
For $q=1$ the theorem specializies to a result of Kleitman \cite{Kle}.

\section{Poincar\'e polynomials for intersections of Schubert varieties}
\noindent

Let $G$ be a complex semisimple algebraic group and $P$ a maximal parabolic subgroup.
In the homogeneous space $G/P$, known as a {\em generalized Grassmannian},
 there are {\em Schubert varieties} $X_u$, indexed by elements $u\in W^J$.
 Here $(W,S)$ is the Weyl group of $G$, $J=S\setminus \{s_{\al}\} $,
  and $W^J$
 the Bruhat poset of minimal coset representatives modulo the corresponding 
 maximal parabolic subgroup $W_J$  of $W$, obtained by removing the node 
 corresponding to $\al$ from
 the Coxeter-Dynkin diagram of $G$.
 
 We refer to the literature, e.g. \cite{BL} and  \cite{Hil}, for definitions and explanations
 of these and other notions concerning algebraic groups, root systems and
 Schubert varieties. For instance, we consider the notion of {\em minuscule weight}
 known, see \cite[p. 21]{BL} or  \cite[pp. 172--174]{Hil}. The Grassmannian $G/P$ is
 {\em of minuscule type} if the fundamental weight $\al$ is minuscule.
 The Grassmannians of minuscule type have been classified. In addition to the
 classical Grassmannians of $k$-dimensional subspaces in $\C^d$
 (the type $A$ case), there are four infinite families and three sporadic cases.

The {\em Poincar\'e polynomial} of a Schubert variety $X_u$ is
$P_{X_{u}}(q)\defeq \sum \be_i(X_u)q^{i}$,
where the Betti numbers are the dimensions of classical singular cohomology
groups with complex coefficients: $ \be_i(X_u) = \dim H^{i}(X_u ; \C)$.
Explicit expressions for  $P_{G/P}(q)$ are known in terms of parameters of
the groups in all cases.

\bthm Let $X_u$ and $X_v$ be Schubert varieties in a Grassmannian
$G/P$ of  minuscule type. Then,
$$ P_{X_u}(q) \cdot P_{X_v}(q)  \ll P_{G/P}(q) \cdot P_{X_u \cap X_v}(q)
$$
\ethm
\bproof
It is a well-known consequence of the decomposition into Schubert cells that 
$\be_{2i}(X_u)$ equals the number of elements 
$w$ such that $\ell(w)=i$ and $w\le u$ in the Bruhat poset $W^J$.
Also, $\be_{2i+1}(X_u)=0$. On the other hand,
it is known from the work of Proctor \cite{Pro} that $W^J$ is a distributive
lattice if and only if $J=S\setminus  \{s_{\al}\}$ for a minuscule weight $\al$.
Being a lattice implies that $X_u \cap X_v =X_{u\wedge v}$.

Hence, the result is obtained by applying  Theorem \ref{qfkg} to $W^J$, letting
$$\mu(w)\equiv 1, \quad\quad 
g(w)=\begin{cases} 
1, \mbox{ if } w\le u \\ 
0, \mbox{ if } w\not\le u \end{cases} , \mbox{ and } \quad
  h(w)=\begin{cases} 1, \mbox{ if } w\le v \\ 0, \mbox{ if } w\not\le v \end{cases}.$$

\eproof

\section{A more general $q$-FKG inequality}


Let $L$ be a finite distributive lattice, and
for $x\in L$  let $m(x)$ denote the number of maximal chains in
the interval $[\hz, x]$.
The following result  is  due to Fishburn \cite{Fis}
in an equivalent version formulated for linear extensions of posets, 
see also \cite[Lemma 3.1]{BT}.


\blem {\rm (Fishburn)} \label{fishburn} 
The function
$x\mapsto \frac{m(x)}{ r(x)!}$ is log-supermodular on $L$. \\  
\elem

We need is a slightly elaborated version.
\bprop \label{fishburn.ext} 
Let $0\le s\le t$ be real numbers.
Then the function
$x\mapsto \frac{m^t(x)}{ (r(x)!)^s}$ is log-supermodular on $L$. \\  
\eprop
\bproof
It follows from the modular law for the rank function (\ref{modular}) that $x\mapsto r(x)!$ is
log-supermodular. Also, log-supermodularity is preserved by powers and products. 
Fishburn's lemma and
$$  \frac{m^t(x)}{ (r(x)!)^s} =  \left(\frac{m(x)}{r(x)!}\right)^t \cdot  (r(x)!)^{t-s}
$$
therefore complete the proof.
\eproof



\bthm \label{qfkg.alt} Let $L$ be a finite distributive lattice, 
$\mu: L\rarr \R^+$ a log-supermodular weight function, and $0\le s\le t$.
For functions $k: L\rarr \R^+$ let
\beqas
E(k;q) & \defeq & \sum_{x\in L} k(x) \mu(x) m^{t}(x) \frac{q^{r(x)}}{(r(x)!)^s} \in \R^+ [q].
\eeqas
If $g, h: L\rarr \R^+$ are comonotone functions, then
$$E(g;q)\cdot E(h;q)\   \ll \ E (1;q)\cdot E(gh;q).
$$
For countermonotone functions the inequality is reversed.
\ethm

\bproof 
The
inequality  follows from Theorem \ref{qfkg} and Proposition \ref{fishburn.ext}.
\eproof

\section{Series weighted by Young tableaux}
\noindent

Young's lattice $Y$ consists of all number partitions ordered by inclusion
of their Ferrers diagrams. It is a distributive lattice.
For two partitions 
$$\la = (\la_1 \ge \la_2 \ge \cdots)
\quad \mbox{ and }\quad  \si = (\si_1 \ge \si_2 \ge \cdots)$$ the lattice operations are
$$\la\vee \si = (\max(\la_1, \si_1) \ge \max(\la_2, \si_2) \ge \cdots)$$
and
$$\la\wedge \si = (\min(\la_1, \si_1) \ge \min(\la_2, \si_2) \ge \cdots).$$
For notation and definitions concerning number partitions, Young tableaux and Young's
lattice we refer to \cite[Ch. 7]{EC2}. 
In particular, $f_{\la}$ denotes the number of standard
Young tableaux of partition  shape $\la$, and  $|\la|=\sum \la_i$.

\bprop\label{young}
Suppose that $0\le s\le t$. Then the function
$$  \la \mapsto \frac{f_{\la}^t}{(|\la| !)^s}$$ 
from $Y$ to $\R^+$ is log-supermodular.
\eprop


Proposition \ref{young} is the Young's lattice special case of Proposition \ref{fishburn.ext} 
and therefore relies on Lemma \ref{fishburn}.
Since the  proofs in  \cite{Fis} and \cite{BT} for Lemma \ref{fishburn}  are somewhat involved,
we provide a simpler direct proof for the Young's lattice special case, 
thus  making this paper self-contained concerning the results in this section.
Our proof is based on the hook formula and Lemma \ref{easy}.

Let $L$ be a distributive lattice, and for $x,y\in L$ let $d(x,y)
\defeq r(x\vee y) - r(x\wedge y)$. Then $d$ is a metric on $L$ (we don't use this fact),
and  $d(x,y)=2$ if and only if $x$ and $y$ are either the top and bottom elements or 
the two middle elements of an interval of length $2$.

\blem \label{easy}
A  function $\mu: L\rarr \R^+ \setminus \{0\}$ is   {\em log-supermodular} if and only if
$$\mu(x) \mu(y) \le \mu(x\wedge y)\mu(x\vee y)   \mbox{ for all } x,y\in L
\mbox{ such that $d(x,y)=2$}. $$
\elem

\bproof
Suppose that the inequality holds for  all $x,y$ such that
$d(x,y)=2$. We prove that it then holds for all pairs $x,y$, by induction
on $d(x,y)$. 

Suppose that $d(x,y)\ge 3$ and (without loss of generality) that
$r(x)-r(x\wedge y)\ge 2$. Let $a$ be an atom in the interval $[x\wedge y, x]$.
Notice that
$d(a,y)<d(x,y)$ and $d(x,a\vee y)<d(x,y)$.
Hence, the induction assumption gives 
$$\frac{\mu(y)}{\mu(x\wedge y)}=
\frac{\mu(y)}{\mu(a\wedge y)} \le
\frac{\mu(a\vee y)}{\mu(a)}\le
\frac{\mu(x\vee y)}{\mu(x)}
$$

\eproof

\noindent
{\em Proof of the Young's lattice case of Lemma  \ref{fishburn}.}



For two partitions $\la$ and $\si$, let
$\si\setminus\la$ denote  the set of cells in the Ferrers diagram
of $\si$ that are not in $\la$. 
By Lemma \ref{easy} we may in the following assume that $|\, \si\setminus\la \, |=1$
and that $|\, \la\setminus\si \, |=1.$
Let 
$s\defeq |\si| -1 =|\la|-1$.

Let $c$ be the cell of the Ferrers diagram  defined by $\si\setminus\la = \{c\}$. Then,
$$\la\wedge \si=\si\setminus c \mbox{ \; and \; } \la\vee\si =\la\cup c.$$
Using the hook formula for $f_{\la}$ (see \cite[p. xx]{EC2}), 
the inequality we have to prove is the following
$$ \frac{1}{\prod_{a\in \la} h_a} \,\cdot\, \frac{1}{\prod_{b\in \si} h_b} \,\le\,
\frac{1}{\prod_{a\in \la\cup c} h_a} \;\cdot\;   \frac{1}{\prod_{b\in \si\setminus c} h_b},
$$
which we rewrite
\beq \prod_{a\in \la} h_a^{\la} \prod_{b\in \si} h_b^{\si} \label{eq1}
\;\ge\; \prod_{a\in \la\cup c} h_a^{\la\cup c} \prod_{b\in \si\setminus c} h_b^{\si\setminus c},
\eeq
where the upper indices of the hook numbers refer to the respective shapes in question.

Let $Z$ be the set of all cells of the Ferrers diagram,
other than $c$,  whose place is in the same row or same column as 
the critical cell $c$. We observe that 
$$ h_a^{\la} =  h_a^{\la\cup c} \mbox{ and } h_b^{\si} = h_b^{\si\setminus c}
\mbox{ for all cells  $a, b$ {\em not} in $Z$,}$$
and 
$$h_a^{\la\cup c}= h_a^{\la}+1  \mbox{ and } h_b^{\si\setminus c}=h_b^{\si} -1 
\mbox{ for all cells  $a, b\in Z$.}$$
Thus, cancelling equal factors in inequality (\ref{eq1}) all that remains to prove is that
\beqs \prod_{a\in Z} h_a^{{\la}}  \prod_{b\in Z} h_b^{\si}
  \;\ge\; \prod_{a\in Z} (h_a^{\la}+1)  \prod_{b\in Z} (h_b^{\si}-1).
\eeqs
This  is 
clear since $\si\setminus c \sbseq \la$ implies that $h_a^{\si}\le h_a^{\la}$, and hence
\beqs h_a^{{\la}}  h_a^{\si}
  \;\ge\; (h_a^{\la}+1)   (h_a^{\si}-1), \mbox{ for all } a\in Z.
\eeqs
\hfill $\Box$

For formal power series $G(z), H(z) \in \R[[z]]$,
$$\mbox{ let $G(z) \ll H(z)$
mean that $H(z)-G(z)\in \R^+ [[z]]$.}$$
We obtain the following correlation-type inequality for monotone functions
on Young's lattice.

\bthm \label{SYTseries}
Suppose that $\mu: Y\rarr \R^+$ is log-supermodular and $0\le s\le t$. For functions
$k: Y \rarr \R^+$ let
$$F(k ; z)\defeq\sum_{\la\in Y} k(\la) \, \mu(\la)\,  (f_{\la})^t \, \frac{z^{|\la|}}{(|\la| !)^s} \, \in \R^+ [[z]].$$
If $g, h: Y \rarr \R^+$ are  comonotone 
functions,  then 
$$F(g;z) \, \cdot \,  F(h;z) \, \ll \, F(1;z)  \, \cdot \, F(gh;z).
$$
For countermonotone functions the inequality is reversed.
\ethm
\bproof
To prove this coefficient-wise inequality for the degree $d$ term it suffices to consider
the truncations of the four power series to polynomials of degree  $d$,
obtained by omitting all higher-degree terms.
Consider the partition $d^d \defeq (d,d,\ldots ,d)$, and let $L\defeq [ \emptyset, d^d]$. 
Then $L$ is an interval in $Y$ and hence a finite distributive lattice.
Now observe that the degree $d$ truncation of the polynomial
$E(g;z)$ of Theorem \ref{qfkg.alt}
computed on $L$ equals  the degree $d$ truncation of  $F(g;z)$,
and similarly for the other four factors. Hence, the result is a consequence
of Theorem \ref{qfkg.alt}.
\eproof

A special case deserves mention.
The function
$$  \la \mapsto \frac{f_{\la}^2}{|\la| !}$$ 
from $Y$ to $\R^+$ is known as {\em Plancherel measure on partitions}, and
$$  \la \mapsto \pi_{\theta}(\la) \defeq e^{-\theta}\frac{\theta^{|\la|} f_{\la}^2}{(|\la| !)^2}$$ 
as its {\em poissonization}, see \cite{BOO}. Here $\theta >0$ is a real parameter.

\bcor
For functions $k: Y \rarr \R^+$ let
$$P(k ; z)\defeq  \sum_{\la\in Y} k(\la) \, \pi_{\theta}(\la)z^{|\la|}, $$
where $\pi_{\theta}$ is the poissonization of Plancherel measure.
If $g, h: Y \rarr \R^+$ are  comonotone 
functions,  then 
\beq\label{corr1} P(g;z) \, \cdot \,  P(h;z) \, \ll \, P(1;z)  \, \cdot \, P(gh;z).
\eeq
  In particular, under condition of convergence there is the correlation  inequality
\beq\label{corr2} \sum_{\la\in Y} g(\la) \, \pi_{\theta}(\la) \, 
 \, \cdot \, 
 \sum_{\la\in Y} h(\la) \, \pi_{\theta}(\la) \,
 \,\le \,  \, \sum_{\la\in Y} g(\la) h(\la) \, \pi_{\theta}(\la) \,
  \eeq
\ecor

We end with a curiosity.
Let $\la\vdash n$ denote that $\la$ is a partition of the integer $n$, 
let $p(n)$ be the number of partitions of  $n$, and let $i(n)$ be the
number of involutions of an $n$-element set.
It is known that
\beqas
\sum_{\la\vdash n} f_{\la}^2 = n! , \qquad
\sum_{\la\vdash n} f_{\la} = i(n),  \qquad
\sum_{n\ge 0} i(n)
\frac{z^n}{n!} = e^{z+z^2 /2} ,
\eeqas
see \cite{EC2}. Thus we have:
$$\sum_{\la\in Y}\,  f_{\la}^{k}\, \frac{z^{|\la|}}{|\la| !}=
\bca
 \sum_{n\ge 0} p(n) \frac{z^n}{n!}, \mbox{ if } k=0  \\
 \\
 e^{z+z^2 /2}, \mbox{ if } k=1 \\ \\
\frac{1}{1-z}, \mbox{ if } k=2.
\eca
$$
Such explicit expressions do not seem to be known for other values of the power $k$.
However, the following relations between  such series can be observed.

The function $\la \mapsto f_{\la}^s$ is decreasing if $s<0$, and
increasing if $s>0$.
Hence, if $st>0$ 
 we obtain this special case of Theorem \ref{SYTseries}:
\beq\label{sample2}
\( \sum_{\la\in Y}  \frac{f_{\la}^{s+1} z^{|\la |}}{|\la| !} \)
\( \sum_{\la\in Y}  \frac{f_{\la}^{t+1} z^{|\la|}}{|\la| !} \) \ll
e^{z +{z}^2/2} \, \cdot \,
 \sum_{\la\in Y}  \frac{f_{\la}^{s+t+1} z^{|\la|}}{|\la| !} ,
\eeq
whereas if $st<0$ the inequality goes the other way.


\section{Remarks}

{\bf 1.} The original FKG inequality has been generalized in several directions. For instance,  Sahi
\cite{Sah} gives a polynomial version different from the one considered here.

Ahlswede-Daykin's  {\em four function theorem} 
\cite{ADa} states the following.
Suppose that $L$ is a finite distributive lattice and that $\al, \be, \ga, \de$ are nonnegative
real functions on $L$ such that
\beq\label{AD}
\al(x)\ \be(y) \le \ga(x\vee y)\ \de(x\wedge y),
\eeq
for all $x,y \in L$. Then
\beqs
\sum_{x\in \AA} \al(x) \ \cdot \ 
\sum_{x\in \BB} \be(x) \ \le \
\sum_{x\in \AA\vee \BB} \ga(x) \ \cdot \ 
\sum_{x\in \AA\wedge \BB} \de(x),
\eeqs
for all subsets $\AA, \BB \sbseq L$. Here,
$\AA\vee \BB \defeq \{x\vee y \st x\in \AA, y\in \BB\}$ and
$\AA\wedge \BB \defeq \{x\wedge y \st x\in \AA, y\in \BB\}$. The FKG
inequality is obtained from this by putting
$$\al=g\mu, \ \
\be=h\mu, \ \
\ga=gh\mu, \ \
\de=\mu \ \mbox{ and $\AA=\BB=L$.}
$$
This prompts the following.

{\bf Question:} Let $\al, \be, \ga, \de: L \rarr \R^+$ be functions as above, satisfying
condition (\ref{AD}). Does it then follow that
\beqs
\sum_{x\in \AA} \al(x) q^{r(x)} \ \cdot \ 
\sum_{x\in \BB} \be(x) q^{r(x)} \ll
\sum_{x\in \AA\vee \BB} \ga(x) q^{r(x)} \ \cdot \ 
\sum_{x\in \AA\wedge \BB} \de(x) q^{r(x)} \ ?
\eeqs

{\bf 2.} We are grateful to  P. Br\"and\'en for bringing references
\cite{BT} and  \cite{Fis} to our attention.


\begin{thebibliography}{}

\bibitem{ADa} R.~Ahlswede and D.~E.~Daykin, An inequality for the weights of
two families of sets, their unions and intersections, {\em Z. Wahrs. Verw. Gebiete}
{\bf 43} (1978), 183--185.

\bibitem{AS} N.~Alon and J.~H.~ Spencer, {\em The Probabilistic Method,
3rd Ed.}, Wiley-Interscience, New York, 2008.

\bibitem{BL} S.~Billey and V.~Lakshmibai, {\em Singular Loci of Schubert Varieties},
Progress in Math. No. 182, Birkh\"auser, Boston, 2000.

\bibitem{Bir} G. Birkhoff, {\em Lattice Theory, 3rd ed.}, Amer. Math. Soc. Colloquium Publ. No. 25,
American Mathematical Society, Providence, RI, 1967.

\bibitem{BOO} A.~Borodin, A.~Okounkov and G.~Olshanski, Asymptotics of Plancherel measures for
symmetric groups, {\em J. Amer. Math. Soc.} {\bf 13} (2000), 481--515.

\bibitem{BT} G.~R.~Brightwell and W.~T.~Trotter, A combinatorial approach to
correlation inequalities, {\em Discrete Math.} {\bf 257} (2002), 311--327.

\bibitem{Fis} P.~C.~Fishburn, A correlational inequality for linear extensions of 
a poset, {\em Order} {\bf 1} (1984), 127--137.

\bibitem{FKG} C.~M.~Fortuin, P.~W.~Kasteleyn, and J.~Ginibre, Corrrelation inequalities on 
some partially ordered sets, {\em Commun. Math. Phys.} {\bf 22} (1971), 89--103.

\bibitem{Hil}  H. Hiller, {\em Geometry of Coxeter Groups}, Research Notes in Math.  No. 54,
Pitman, Boston, 1982.

\bibitem{Kle} D.~J.~ Kleitman, Families of non-disjoint subsets, {\em J. Combinat. Theory}
{\bf 1} (1966), 153--155.

\bibitem{Pro} R.~{A}.~ Proctor, Bruhat lattices, plane partition generating functions, and
minuscule representations, {\em Europ. J. Combinatorics} {\bf  5} (1984), 331--350. 

\bibitem{Sah} S. Sahi, The FKG inequality for partially ordered algebras, {\em J. Theor. Probab.}
{\bf 21} (2008), 449--458

\bibitem{EC1} R.~P.~Stanley, {\em Enumerative Combinatorics, Vol.  1}, Cambridge Studies
in Advanced Mathemtics No. 49, Cambridge University Press, Cambridge, UK, 1997.

\bibitem{EC2} R.~P.~Stanley, {\em Enumerative Combinatorics, Vol.  2}, Cambridge Studies
in Advanced Mathemtics No. 62, Cambridge University Press, Cambridge, UK, 1999.


\end{thebibliography}
\end{document}